\documentclass[12pt,russian]{article}

\usepackage[cp1251]{inputenc}
\usepackage[english,russian]{ babel}
\usepackage[centertags]{amsmath}
\usepackage{amsfonts}
\usepackage{amssymb}
\usepackage{amsthm}
\usepackage{amsfonts,amsmath,amsthm}
\usepackage{latexsym,verbatim}

\theoremstyle{definition}

\theoremstyle{plain}
\newtheorem{thm}{Theorem}
\newtheorem{cor}{Corollary}
\newtheorem{lem}{Lemma}

\newcommand{\eps} {\varepsilon}
\newcommand{\R} {\mathbb{R}}

\newcommand{\CC} {\mathbb{C}}

\newcommand{\Ga} {\widehat{\Gamma}}
\def\d{\delta}
\def\r2{\root\of{2}}
\def\sp{{\rm sp}}
\def\arg{{\rm arg}}
\def\const{{\rm const}}
\def\grad{{\rm grad}}
\def\Int{{\rm Int}}
 \def\re{{\rm Re}}
 \def\im{{\rm Im}}

\title{
\begin{center}
A multidimensional version of Levin's Secular Constant Theorem and
its applications.
\end{center}}
\author{S.Yu. Favorov,
 N. Girya
}

\date{}
\begin{document}
\maketitle
\begin{abstract}
We study holomorphic almost periodic functions on a tube domain
with the spectrum in a cone. We extend to this case Levin's
theorem on a connection between the Jessen function, secular
constant, and the Phragmen-Lindel\"of indicator. Then we obtain a
multidimensional version of  Picard's theorem on exceptional
values for our class.
\end{abstract}

An almost periodic function with  bounded from below spectrum has
some specific properties. Namely, it  extends to the upper
half-plane  as a holomorphic almost periodic function $f$ of
exponential type (H.Bohr \cite{bor}), then  $\log|f|$ and the mean
value of $\log|f|$ over a horizontal line (so-called Jessen's
function) are  of the same growth along the imaginary positive
semi-axis (B.Jessen, H.Tornehave \cite{JessTorn} and B.Ja.Levin
\cite{levin}). The last result (together with the discovered by
Ph.Hartman \cite{Hatr}, and B.Jessen, H.Tornehave \cite{JessTorn}
connection between Jessen's function, mean motions of $\arg f(z)$,
and distribution of zeros for holomorphic almost periodic
functions on a strip) shows the regularity of functions from that
important class.

 In the end of the last century,  L.I.Ronkin \cite{ron},
\cite{ron92}, \cite{R} created the theory of holomorphic almost
periodic functions and mappings defined on tube domains of the
multidimensional complex space. Introduced by him Jessen's
function of  several variables plays the main role in  value
distribution theory  for almost periodic holomorphic mappings.

Here we continue the line of investigation in \cite{FaUd} and
\cite{gifa} of the class of almost periodic functions on a tube
domain with the spectrum in a cone. Namely,  we find a connection
between asymptotic behavior of  Jessen's function and the polar
indicator. Then we introduce a multidimensional analogue of the
secular constant and study its asymptotic behavior. Also, we
obtain  a multidimensional version of  Picard's Theorem on
exceptional values for our class.

\bigskip

 Let us give a more detailed description of the subject.

Suppose $f$ is a $2\pi$-periodic function with the convergent
Fourier series
  $f(x)= \sum\limits_{n\ge n_0}a_n e^{inx},\,n_0\le0,\,a_{n_0}\neq 0$.
Then $f(z)= \sum\limits_{n\ge n_0}a_n e^{inz}$, $z=x+iy,$ is  a
natural extension of $f(x)$ to the upper half-plane
$\mathbb{C^+}$. Clearly, $f(z)$ is  a holomorphic function  of
exponential type $|n_0|$ without zeros in some half-plane $y>y_0$
and
 $$
\lim_{y\to+\infty}y^{-1}\frac{1}{2\pi}\int_{-\pi}^\pi
\log|f(x+iy)|dx= \lim_{y\to+\infty}y^{-1}\log|f(iy)|  =-n_0.
 $$

In \cite {bor} and \cite{JessTorn},  these properties were
generalized to almost periodic functions $f$ with  bounded from
below spectrum under the condition $\Lambda^0=\inf\sp f\in\sp f$.
 One should only replace  the mean value over the period by
{\it Jessen's} function
 \begin{equation}\label{Jessen}
J_f(y)=\lim_{S\to\infty}(2S)^{-1}\int_{-S}^S \log|f(x+iy)|dx;
 \end{equation}
 the number $n_0$  by $\Lambda^0$, and make use of the
Phragmen-Lindel\"of Principle (see  a footnote in the proof of
Theorem 1).

Note that the limit in (\ref{Jessen}) exists for every holomorphic
almost periodic function on a strip $\{z=x+iy:\,a<y<b\}$ and the
function $J_f(y)$ is convex on $(a,\,b)$. Then for all
$y\in(a,\,b)$,  maybe with except of some countable set $E_f$, we
have
 \begin{equation}\label{mot}
  J_f'(y)=-c_f(y),
\end{equation}
 where
 $$
c_f(y)=\lim_{\gamma-\beta\to\infty}\frac {\arg f(\gamma+iy)-\arg
f(\beta+iy)}{\gamma-\beta}
 $$
 is the {\it mean motion}, or {\it secular number}, of the function $f$;
here $\arg f(x+iy)$ is a continuous branch of the argument of $f$
on the line $y=\const$. By the way, equality (\ref{mot}) and the
Argument principle imply that the number $N(-S,\,S,\,y_1,\,y_2)$
of zeros of the function $f$ in the rectangle
$\{|x|<S,\,y_1<y<y_2\}$ \footnote{Zeros should be counted with
multiplicities}  has a density
\begin{equation}\label{dense}
\lim_{S\to\infty}(2S)^{-1}N(-S,\,S,\,y_1,\,y_2)=J_f'(y_2)-J_f'(y_1)
\end{equation}
for all $y_1,\,y_2\not\in E_f$. It can  also be proved that $f$
has no zeros on a substrip $\{\alpha<y<\beta\}$ if and only if
$J_f(y)$ is a linear function on the interval $(\alpha,\,\beta)$.
In this case,
$$
f(z)=e^{ic_fz+g(z)},
$$
where $g(z)$  is almost periodic on the strip
$\{z=x+iy:\,x\in\R,\,\alpha<y<\beta\}$.

 Thus, an almost periodic function $f$
with the property $-\infty<\Lambda^0=\inf\sp f\in\sp f$ is
extended to $\CC^+$ as  a holomorphic almost periodic  function.
Then we get
\begin{equation}\label{main}
-\Lambda^0=\lim_{y\to+\infty}\frac{\log|f(iy)|}{y}=
\lim_{y\to+\infty}\frac{J_f(y)}{y}=
\lim_{y\to+\infty}J_f'(y)=-\lim_{y\to+\infty}c_f(y)
\end{equation}
(see, for example, \cite {JessTorn}, \cite{lev}).

In the case $\Lambda^0\not\in\sp f$,  the function is also
extended to  $\CC^+$ as  a holomorphic almost periodic function;
the equalities (\ref{main}) are also valid, but the proof of the
second equality is complicated,  and  this is the contents of
Levin's Secular Constant Theorem \cite{levin}, \cite{lev}.

\medskip

Note that there exists  a natural connection between  the
distribution of zeros of an almost periodic holomorphic function
on the upper half-plane and  the configuration of its spectrum:

{\bf Theorem B} (\cite{b}). {\it Suppose that the spectrum $\sp f$
of an almost periodic function $f$ on $\CC^+$ is bounded from
below. Then
\begin{enumerate}
\item if $\Lambda^0=\inf\sp f\ge 0$, then $f(z)$ tends to a finite
limit as  $y \to \infty$ on $\CC^+$ uniformly in $x\in\R$,

 \item if $\Lambda^0=\inf\sp f<0$ and $\Lambda^0 \in\sp f$,
then $f(z)\to\infty$ as  $y \to \infty$ on $\CC^+$ uniformly in
$x\in\R$,

 \item if $\Lambda^0=\inf\sp f<0$ and $\Lambda^0\not\in\sp f$,
then the function $f(z)$ takes every complex  value on the
half-plane $y>q\ge 0$ for each $q < \infty$.
 \end{enumerate}}

\bigskip

 To discuss the multidimensional case,  we need the following
definitions.

 Let $z=(z_1,\dots,z_n)\in\mathbb{C}^p$,
$z=x+iy\in\mathbb{C}^p$, $x\in\R^p,\,y\in\R^p$. By  $\langle x,\,y
\rangle$ or $\langle z,\,w
 \rangle$ denote
the scalar product (or the Hermitian scalar product
  for $z,w \in\mathbb{C}^p$).
By $|.|$  denote the Euclidean norm  on $\R^p$ or $\mathbb{C}^p$.
Also, for $x=(x_1,x_2,\dots,x_p)$ put $'x=(x_2,\dots,x_p)$.
Further, by $T_K$  denote a tube  set $$
T_K=\{z=x+iy\in\mathbb{C}^p:\, x\in \R^p, y\in K \},
$$
where $K\subset\R^p$  is the base of the tube  set.

A vector $\tau\in\mathbb{R}^p$ is called an $\eps$-almost period
of a function $f(z)$ on $T_K$ if
$$
  \sup_{z\in T_K} |f(z+\tau) -f(z)| < \eps.
$$
The function $f$ is called {\it almost periodic on} $T_K$ if for
every $\eps>0$ there exists $L=L(\eps)$ such that every
$p$-dimensional cube in $\mathbb{R}^p$ with the side of  length
$L$ contains at least one $\eps$-almost period of $f$. In
particular, when $K=\{0\}$, we get the definition of an almost
periodic function on $\mathbb{R}^p$.

A function $f(z)$, $ z\in T_\Omega$,
 where $\Omega$  is a domain in $\mathbb{R}^p$,
is called almost periodic if its restriction to $T_K$ is  an
almost periodic function for every compact set $K\subset\Omega$.

The spectrum $\sp f$ of an almost periodic function $f(z)$ on
$T_K$ is the set of vectors $\lambda \in \R^p$ such that the
Fourier coefficient
\begin{equation}\label{sr}
a_{\lambda}(y, f)=\lim_{S\to\infty} \frac{1}{(2S)^p}\int\limits
_{|x_j|<S, j=1..p} f(x+iy) e^{-i \langle x, \lambda \rangle}
dm_p(x)
\end{equation}
does not vanish on $K$; here $m_p$ is the Lebesgue measure on
$\R^p$. The spectrum of every almost periodic function $f$ is at
most countable, therefore we have
 $$
f(x+iy) \sim \sum a_n(y) e^{i \langle x, \lambda^{n} \rangle},
 $$
where $\{\lambda^{n}\}_{n\in \mathbb{N}}=\sp f$ and
$a_n(y)=a_{\lambda^{n}}(y, f)$. Note that for any  given countable
set $\{\lambda^{n}\}$
 the function
$\sum_{n\in\mathbb{N}}n^{-2}e^{i\langle x,\lambda^{n}\rangle}$ is
almost periodic on $\R^p$ with the spectrum $\{\lambda^{n}\}$.

 In \cite{ron} L.I. Ronkin introduced the notion of Jessen's function
of an almost periodic holomorphic function  $f$ on
 $T_{\Omega}$ by  the formula
 $$
J_f(y) = \lim\limits_{S \to \infty} \frac 1 {(2S)^p}
\int_{[-S,S]^p}\log|f(x+iy)|dm_p(t).
 $$
 Using the methods of the theory of distributions and special properties
of zero sets for holomorphic functions, L.I.Ronkin checked that
the limit exists and defines a convex function in $y\in\Omega$. He
also established the multidimensional analogue of equality
(\ref{dense})
 $$
\lim\limits_{S\to\infty}{m_{2p-2}\{z=x+iy:\,x\in
[-S,\,S]^p,\,y\in\omega,\,f(z)=0\}\over(2S)^p}=\kappa_p\mu_J(\omega),
$$
where $\mu_J$ is the Riesz measure of $J(y)$,
$\omega\subset\overline{\omega}\subset\Omega$,
$\mu_J(\partial\omega)=0$, and the area of the zero sets  is taken
counting the multiplicity.

Also, in \cite{ron1} L.I.Ronkin proved that the products
$b_n(y)=a_n(y)e^{ \langle y, \lambda^n \rangle}$ do not depend on
$y$ for every holomorphic almost periodic function $f(z)$ on
$T_\Omega$; in particular, the  coefficient $b_0$ corresponding to
the exponent $\lambda=0$  does not depend on $y$. In the case, the
Fourier series turns into the Dirichlet series
 \begin{equation}\label{Dir}
f(z) \sim \sum_{\lambda^{n} \in \R^p} b_n e^{i \langle z,
\lambda^{n} \rangle},\quad b_n\in\mathbb{C}.
\end{equation}

In \cite{ron92} L.I.Ronkin obtained the following results.

{\bf Theorem R.} {\it Let $f$ be a holomorphic almost periodic
function on $T_{\Omega}$. Then the function $J_f(y)$ is linear on
the domain $\Omega'\subset\Omega$ if and only if the function $f$
has no zeros in $T_{\Omega'}$. Moreover, in this case
 \begin{equation}\label{sec}
  f(z)=\exp\{i\langle c_f,\,z\rangle+g(z)\},\quad z\in
T_{\Omega'},
\end{equation}
where $c_f\in\R^p$ and $g(z)$  is an almost periodic function on
$T_{\Omega'}$.}

In conditions of Theorem R,  we have
$$
J_f(y)=-\langle c_f,\,y\rangle+\re\, b_0,\quad y\in\Omega',
$$
where $b_0$ is the corresponding  coefficient of the
Dirichlet-series expansion of the function $g$. Therefore, the
following definition seems to be natural:

\medskip
{\bf Definition.} {\it The function $-\grad J_f(y),\, y\in\Omega,$
is the secular vector of the almost periodic holomorphic function
$f$ on $T_\Omega$}.
\medskip

In order to formulate our results,  we need some definitions and
notations.

 A cone $\Gamma\subset\mathbb{R}^p$ is the set with the property
 $ y\in\Gamma, t>0
\Rightarrow ty\in\Gamma. $ We will consider convex cones with
non-empty interior and  such that $\overline{\Gamma}\bigcap
(-\overline{\Gamma})=\{0\}$. By $\Ga$ denote the conjugate cone to
$\Gamma$, i.e., $\Ga = \{ x \in \R^p : \langle x, y \rangle \ge
0\;\;\; \forall y\in~\Gamma\}$;  note that
$\widehat{\Ga}=\overline\Gamma$. As usual, $\Int A$ is the
interior of the set $A$,  and $H_E(x)=\sup_{\lambda\in E}\langle
x,\,\lambda \rangle$ is the {\it support function} of the set
$E\subset\R^p$.

  Let $f$ be a holomorphic almost periodic function on a tube $T_\Gamma$
with  an open cone $\Gamma$ in the base. By definition, put
 $$
  h_f(y)=\sup_{x\in\mathbb{R}^p}\overline{\lim_{r\to\infty}}
\frac{\ln|f(x+iry)|}{r},  \quad y \in \Gamma.
 $$
The function $h_f$ is called the {\it $P$-indicator} of $f$ (see
\cite{R}, p.245).

   {\bf Theorem A.} \label{b} (\cite{gifa})
 {\it Let $\Gamma$ be a closed cone in $\R^p$, and $f(x)$ be an
almost periodic function on $\R^p$. Then $f$ is extended
holomorphically to $T_{\Int\Ga}$ with the estimates
\begin{equation}\label{exp_rost} \exists b < \infty \quad \forall
\Gamma'=\overline{\Gamma'} \subset\Int\Ga\cup\{0\} \quad \exists
B(\Gamma') \quad \forall z \in T_{\Gamma'} \quad |f(z)| \le
B(\Gamma') e^{b|y|},
\end{equation}
if and only if $spf\subset\Lambda+\Gamma$ for some
$\Lambda\in\mathbb{R}^p$.  If this is the case, then $f(z)$ is
almost periodic on $T_{\Int\widehat{\Gamma}}$ and for all $y\in
\Int\widehat{\Gamma}$
\begin{equation}\label{exp_rost1}
h_f(y)=H_{\sp f}(-y).
\end{equation}}
 For almost periodic functions with  bounded spectrum,  equality
(\ref{exp_rost1}) was proved in \cite{FaUd}.

 \bigskip

 The following theorem is the main result of our paper.

\begin{thm} \label{thm1}
Let $\Gamma$ be a closed cone in $\R^p$, and $f(x)$ be an almost
periodic function on $\R^p$ such that $f$ is extended
holomorphically to $T_{\Int\Ga}$ with estimates (\ref{exp_rost}).
Then for all $y\in\Int\Ga$
 \begin{equation}\label{th}
 \lim_{R\to\infty} \frac{J_f (R y)} {R} =h_f (y).
 \end{equation}
 Furthermore, the secular vector $-\grad J_f(Ry)$ tends to $\grad
H_{\sp f}(-y)$ as $R\to\infty$ in the sense of distributions.
\end{thm}

{\bf Remark}. Since $J_f(y)$ is  a convex function, we see that
the secular vector is a locally integrable function on $\Int\Ga$.

\medskip
\begin{proof}
 From the beginning  assume that $y^0 = (1,0,0,...0) \in\Int\Ga$, and
  we will prove (\ref{th}) for $y=y^0$.

Put $F(z)=f(z) e^{i\langle
z,\,h_f(y^0)y^0\rangle}\left(\sup_{x\in\R^p}|f(x)|\right)^{-1}$,
$u(z)= \log|F(z)|$. Note that $F(z)$ is an almost periodic
holomorphic function on $T_{\Int\Ga}$ and $|F(x)|\le 1$ on $\R^p$.
  Applying the Phragmen-Lindel\"of Principle
 \footnote{Suppose $g(z)$ is continuous on $\overline{\mathbb{C^+}}$,
holomorphic on $\mathbb{C^+}$, and  bounded on $\R$ function,
which satisfies the condition $\log^+|g(z)|=O(|z|)$ as
$|z|\to\infty$; then for $z=x+iy\in\mathbb{C^+}$ we have
$|g(z)|\le\sup_{x\in\R}|g(x)|e^{\sigma^+y}$, where
 $\sigma^+
=\limsup_{y\to+\infty} y^{-1}\log|g(iy)|$
  (see \cite{koos}, p. 28).}
on the complex one--dimensional  plane
$\{x+wy:\,w\in\mathbb{C^+}\}$, we get
 \begin{equation}\label{Fr-L}
  u(x+ity)\le h_F(y)\,t,\quad\forall\,z=x+iy\in T_{\Int\Ga},\quad t>0.
\end{equation}
Then
  \begin{equation}\label{rost}
h_{F}(y)=h_f(y)-\langle y, h_f(y^0) y^0\rangle,\qquad
h_{F}(y^0)=0.
 \end{equation}
Take  $y=y^0$ in (\ref{Fr-L}). We get
 \begin{equation}\label{ine}
u(z_1,'x) \le 0\quad \forall\quad (x_1,\,'x)\in\R^p,\quad y_1\ge0.
 \end{equation}
Fix $\eps>0$. Since
$\sup_{x\in\R^p}\varlimsup_{r\to\infty}r^{-1}u(x+iry^0)=0$, we see
that for some $x^0=x^0(\eps)\in \R^p$, $r=r(\eps)>0$,
 \begin{equation}\label{3}
u(x^0+iry^0)\ge -\eps r.
 \end{equation}

Using the Poisson formula for the disc
$D(x_1^0+iR,\,R)=\{z_1:\,|z_1-x_1^0-iR|<R\}\subset\CC^+$ with
$R>r$, inequality (\ref{ine}), and Maximum Principle for the
subharmonic function $u(z_1,'x^0)$, we obtain
 $$
u(x_1^0+ir, 'x^0 )\le
 $$
 $$
\le \frac {1}{2\pi} \int _{0}^{2\pi} u(x_1^0+iR+Re^{i\psi},
'x^0)\frac {R^2 - (R-r)^2}{R^2-2R(R-r)cos(\pi/2+\psi)+(R-r)^2}
d\psi
 $$
 $$
\le \frac {r}{4\pi R} \int _{\pi /4}^{3\pi /4}
u(x_1^0+iR+Re^{i\psi}, 'x^0)
 d\psi
\le r(8R)^{-1} \sup_{\psi\in [\pi /4;3\pi /4]}
u(x_1^0+iR+Re^{i\psi}, 'x^0).
$$
Hence (\ref{3}) implies that $u(x_1^0+iR+Re^{i\psi_0},'x^0)\ge
-8\eps R$ for some $\psi_0\in [\pi /4,\,3\pi /4]$. The function
$u(z_1,'x^0)$ is subharmonic in $z_1\in\mathbb{C^+}$. Taking into
account (\ref{ine}) and  the embeddings
 $$
D(x_1^0+2iR,\,R)\subset D(x_1^0+iR+Re^{i\psi_0},\,R+R~/\root\of
{2})\subset\CC^+,
 $$
we get
 $$
-8\eps R\le\frac{2}{\pi
R^2(3+2\,\r2)}\int\limits_{D(x_1^0+iR+Re^{i\psi_0},\,R+R/\root\of{2})}
u (z_1, 'x^0) dm_2(z_1)
 $$
 \begin{equation}\label{4a}
<\frac{1}{3\pi
R^2}\int\limits_{D(x_1^0+2iR,\,R)}u(z_1,'x^0)dm_2(z_1).
 \end{equation}
Remind  that this inequality is valid for all $R>r$.

Put $u_R(z)=R^{-1}u(Rz)$. From (\ref{4a}) it follows that
 \begin{equation}\label{5}
\int\limits_{D(x_1^0/R+2i,\,1)}u_R(z_1,'x^0/R)dm_2(z_1)>
-24\pi\eps.
 \end{equation}

Furthermore,  Theorem A implies that  the function $h_f(y)$ is
continuous. Since (\ref{rost}), we get $h_{F}(y) < \eps$ for
$|y-y^0|<p\delta $ with some
$\delta=\delta(\eps)\in(0,\,1/(p+2))$. If we replace in
(\ref{Fr-L}) $y$ by $y/|y|$, $x$ by $Rx$, and $t$ by $R|y|$, we
obtain
\begin{equation}\label{7}
  u_R(z)=R^{-1} u(Rx+iRy)\le \eps|y|
\end{equation}
for all $z$ from the tube domain
 $$
T^\delta=\{z=x+iy:\,x\in\R^p,\,\left|y/|y|-y^0\right|<p\delta\}.
 $$
By definition, put
 $$
A(x^1)=D(x_1^1+2i,\,1)\times D(x_2^1,\d)\times
D(x_3^1,\delta)\times\dots\times D(x_p^1,\d).
 $$
It can easily be checked that for all $x^1=(x_1^1,'x^1)\in\R^p$ we
have $A(x^1)\subset T^\d$. Also, we may assume that
$\overline{T^\delta}\subset T_{\Int\Ga}\cup\{0\}$. Then for all
$z_1\in\ D(x_1^0/R+2i,\,1)$ the function $u_R(z)$ is subharmonic
in $z_2\in D(x_2^1,\d),\,z_3\in D(x_3^1,\delta),\dots z_p\in
D(x_p^1,\d)$. Hence (\ref{5}) implies
 \begin{equation}\label{10}
\int\limits_{A(x^0/R)}u_R(z)dm_{2p}(z)> -24\d^{2p-2}\pi^p\eps.
 \end{equation}

 Suppose that for some $\tau\in\R^p$ we have
 $$
|F(x^0+\tau+iry^0)-F(x^0+iry^0)|\le e^{-\eps r}-e^{-2\eps r}.
 $$
Then $|F(x^0+\tau+iry^0)|\ge e^{-2\eps r}$ and
$u(x^0+\tau+iry^0)\ge -2\eps r$. Using the  latter inequality
instead of (\ref{3}), we obtain the  relation
 \begin{equation}\label{10a}
\int\limits_{A(x^0/R+\tau/R)}u_R(z)dm_{2p}(z)>
-48\d^{2p-2}\pi^p\eps.
 \end{equation}

Put $u^+_R(z)=\max\{u_R(z), 0\},\,u^-_R(z)=\max\{-u_R(z), 0\}$.
From (\ref{7}) it follows that for all $x^1\in\R^p$ and all $z\in
A(x^1)$ we have
 \begin{equation}\label{9}
u_R(z)<\root\of{10}\eps.
  \end{equation}
Therefore, by (\ref{10a}),
 $$
\int\limits_{A(x^0/R+\tau/R)}u_R^-(z)dm_{2p}(z)= \int\limits_{
A(x^0/R+\tau/R)}u_R^+(z)dm_{2p}(z)-
$$
\begin{equation}\label{11}
\int\limits_{A(x^0/R+\tau/R)}u_R(z)dm_{2p}(z)\le
52\d^{2p-2}\pi^p\eps.
 \end{equation}

In the sequel we need the following lemma:
\begin{lem}
  Let $g(x)$ be an almost periodic function in $x\in\R^p$. Then
for any $\eta>0$ there exist a real $L=L(\eta)$ and a set
$E=E_1\times\dots\times E_p,\,E_j\in\R,$ such that
$E_j\cap[a,\,a+L]\neq\emptyset$  for every $a\in\R,\,j=1,\dots,p$,
and each $\tau\in E$ is an $\eta$-almost period of $g$.
\end{lem}
 \begin{proof} By  Bochner's criterium\footnote{For almost
periodic functions on $\R$ see \cite{lev}, Ch.VI, \S1, or
\cite{C}, p.14-16; the proof for  the multidimensional case is
similar.},  any sequence $t_n\in\R$  has a subsequence $t_{n'}$
such that the functions $g(x+(t_{n'},'0))$  converge uniformly in
$x\in\R^p$. In other words, the functions $g(x_1+t_{n'},'x)$
converge uniformly in $x_1\in\R$ and $'x\in\R^p$. By Bochner's
criterium, the function $g(x_1,'x)$ is almost periodic in
$x_1\in\R$ uniformly in $'x\in\R^{p-1}$. Hence there exist
$E_1\in\R$ and $L=L(\eta)$ such that
$E_1\cap[a,\,a+L]\neq\emptyset$  for all $a\in\R$ and
 $$
 |g(x_1+t,'x)-g(x,'x)|<\eta/p\quad\forall
x_1\in\R,\quad\forall'x\in\R^{p-1},\quad\forall t\in E_1,
 $$
i.e., each $\tau=(t,'0)$ for $t\in E_1$ is an $\eta/p$-almost
period of $g(x)$. In the same way, we find $E_2,\dots,E_p$. It is
clear that every point of $E_1\times\dots\times E_p$ is an
$\eta$-almost period of $g$.
 \end{proof}

Take $S<\infty$, and let $L=L(\eps,\,r)$ be the real from the
Lemma 1. It is not hard to prove that if $R>L~\r2$, then there
exist $\tau^1_1,\dots,\tau^{N_1}_1\in E_1$, $N_1\le 2~\r2 S+2$,
such that
\begin{equation}\label{tau1}
\bigcup_{m=1}^{N_1}\left(\frac{x_1^0+\tau_1^m}{R}-\frac{\r2}{2},
\frac{x_1^0+\tau_1^m}{R}+\frac{\r2}{2}\right)\supset [-S,\,S],
\end{equation}
and each point of $[-S,\,S]$ is contained  in at most two
intervals.
 For the same reasons, if $R>L~\r2/\d$, then for
$j=2,\dots,p$ there exist $\tau^1_j,\dots,\tau^{N_j}_j\in E_j$,
$N_j\le(2~\r2 S+2)/\d$, such that
\begin{equation}\label{tau2}
\bigcup_{m=1}^{N_j}\left(\frac{x_j^0+\tau_j^m}{R}-\frac{\d~\r2}{2},
\frac{x_j^0+\tau_j^m}{R}+\frac{\d~\r2}{2}\right)\supset [-S,\,S].
\end{equation}

Let $F=\{\tau=(\tau_1^{m_1},\dots,\tau_p^{m_p}):\, 1\le m_1\le
N_1,\dots,1\le m_p\le N_p\}$. Note that $F$ contains at most
$(2~\r2 S+2)^p\d^{1-p}$ elements. By definition, put
 $$
\Pi(S,\,\d)=\left\{x+iy:\,x\in[-S,\,S]^p,\,
|y_1-2|<\frac{1}{\r2},\,|y_j|<\frac{\d}{\r2},\,j=2,\dots,p\right\}
 $$
Combining (\ref{tau1}) and (\ref{tau2}), we get
\begin{equation}\label{tau}
 \bigcup_{\tau\in F} A\left(\frac{x^0+\tau}{R}\right)\supset\Pi(S,\,\d).
\end{equation}
Applying  Lemma 1 to the function $F(x+iry^0)$ with $\eta=e^{-\eps
r}-e^{-2\eps r}$ and using (\ref{11}) for every $\tau\in F$, we
obtain
 $$
\int\limits_{\Pi(S,\,\d)}u^-_R(z)dm_{2p}(z)\le \sum_{\tau\in
F}\int\limits_{A(x^0+\tau)}u^-_R(z)dm_{2p}(z)\le 52(2~\r2
S+2)^p\d^{p-1}\pi^p\eps.
$$
Therefore, we have
 \begin{equation}\label{J}
\overline{\lim_{S\to\infty}}\frac 1 {(2S)^p} \int_{\Pi(S,\,\d)}
u_R(z)dm_{2p}(z)\ge -52(\r2\pi)^p\d^{p-1}\eps.
\end{equation}
 It follows from the definition of Jessen's function that
 \begin{equation}\label{jes}
\overline{\lim_{S\to\infty}}\frac 1 {(2S)^p}
\int\limits_{[-S,\,S]} u_R(x+iy)dm_{p}(x)={J_F(Ry)\over R}.
\end{equation}
 The functions $u_R(z)$ are uniformly bounded from above for
$z\in T_{\delta}$. Applying the Fatou lemma to inequality
(\ref{J}), we get
\begin{equation}\label{15}
\int\limits_{|y_1-2|<\frac{\sqrt{2} }{2},|y_2|<
\frac{\delta}{\r2},...|y_p|< \frac{\delta}{\r2}} J_{F} (Ry) d y
\ge -52(\r2\pi)^p\d^{p-1}\eps R,
\end{equation}
for all $R>R(L,\,\d,\,r,\,\eps)$.

To finish the proof,  we need the following simple lemma.
\begin{lem}\label{lem1}
Let $g(t)$  be a convex negative function on $[-\alpha,\,\alpha]$.
Then $g(0)\ge \alpha^{-1} \int_{-\alpha}^\alpha g(t) dt.$
\end{lem}
\begin{proof}
The assertion of the Lemma follows immediately from the inequality
 $$
g(t)\le g(0)\min\{1-t/\alpha,\,1+t/\alpha\}.
 $$
 \end{proof}
Note that (\ref{9}) and (\ref{jes})  imply
\begin{equation}\label{ineq}
J_{F}(R y) \le \root\of{10}\eps R
\end{equation}
for all $y=(y_1,\dots,y_p),\,|y_1-2|<1,\,|y_j|<\d,\,j=2,\dots,p.$
Further, Jessen's function $J_{F}(R y)$ is convex in $y$
(\cite{ron}). Therefore the  function
 $$ g('y) =
\int\limits_{|y_1-2|<\frac{1}{\r2}}J_{F}(Ry)dy_1 - 2\sqrt{5}R\eps
 $$
satisfies the conditions of Lemma 2 in each variable
$y_2,\dots,y_p$ with $\alpha=\d/~\r2$. Applying the  lemma $p-1$
times and using inequality (\ref{15}), we obtain
 $$
\int\limits_{|y_1-2|<\frac{1}{\r2}}
 J_F(Ry_1,'0)dy_1 \ge -40(2\pi)^pR\eps.
$$
Since (\ref{ine}), we see that the integrand is negative.
Moreover, it is convex, therefore $J_F(Ry_1,'0)$ is a
monotonically decreasing function in $y_1$. Then we have
$$
 J_{F}((2-1/\r2)Ry^0)\ge -30(2\pi)^pR\eps.
$$
The  inequality is valid for all $R>R(\eps)$ and $\eps>0$. Thus we
have
\begin{equation}\label{last}
\lim_{R \to \infty} \frac{J_{F}(R y^0)}{R} = 0.
\end{equation}
Since $J_F(y)=J_{f}(y)-\langle y, h_f(y^0)y^0\rangle$, we obtain
(\ref{th}) for $y=y^0$.

For an arbitrary $y'\in\Int\Ga$ consider an orthogonal operator
$A:\R^p\to\R^p$ such  that $A(y^0)=y'$. Put $f_1(z)=f(Az)$. Since
$h_{f_1}(y^0)=h_f(y')$ and $J_{f_1}(y^0)=J_f(y')$, we obtain
(\ref{th}) for $y=y'$.

Further,  from (\ref{Fr-L}) and Theorem A it follows that the
function $J_f(Ry)/R$ is bounded from above on every compact subset
of $\Int\Ga$. Then fix $y^1\in\Int\Ga$ and take $s>0$ such that
$\{y:\,|y-y^1|\le s\}\subset\Int\Ga$. Whenever $|y-y^1|<s$, we
have
$$
2J_f(Ry^1)\le J_f(R(2y^1-y))+J_f(Ry)
$$
and
$$
{J_f(Ry)\over R}\ge 2\inf_{R\ge 1}\left|{J_f(Ry^1)\over
R}\right|-\sup_{R\ge 1}\sup_{|y-y^1|\le
s}{\max\{J_f(Ry),\,0\}\over R}.
$$
 This means that the functions $J_f(Ry)/R$ are uniformly bounded from below on every
compact subset of $\Int\Ga$. Using (\ref{th}) and  the Lebesgue
theorem, we obtain
$$
\int{J_f(Ry)\over R} \varphi(y)dm_p(y)\to\int
h_f(y)\varphi(y)dm_p(y) \quad\hbox{as}\quad R\to\infty
$$
for every test--function $\varphi$ on $\Int\Ga$, i.e., (\ref{th})
is valid in the sense of distributions as well. Therefore,
$$
\grad J_f(Ry)\to\grad h_{f}(y)\quad\hbox{as}\quad R\to\infty
$$
in the sense of distributions and Theorem A implies the last
assertion of Theorem 1.
\end{proof}

 \begin{cor}\label{cor}
Suppose that all conditions of Theorem \ref{thm1} are fulfilled.
If $H_{\sp f}(y)$ is nonlinear on $(-\Ga)$, then $f(z)$
 has zeros on the set $\Int T_{\Ga\cap\{|y|>q\}}$
for each $q<\infty$.
\end{cor}
 \begin{proof}
Theorem A yields that the function $h_f(y)$ is nonlinear for
$y\in\Int\Ga$. Now Theorem \ref{thm1} implies that $J_f(y)$ is
nonlinear on the set $\{Int\Ga\cap\{|y|>q\}\}$ for each
$q<\infty$. Then Theorem R implies that $f(z)$ has zeros on $\Int
T_{\Ga\cap\{|y|>q\}}$.
 \end{proof}

{\bf Applications to distribution of values}.  Here we apply
Theorem 1 to prove the multidimensional variant of Theorem B:

\begin{thm}\label{thm2}
Let $\Gamma \subset \R^p$ be a closed convex cone  and $f(x)$ be
an almost periodic function on $\R^p$  that has a holomorphic
extension $f(z)$ to $T_{\Int\Ga}$ with estimates (\ref{exp_rost}).
Then
 \enumerate
 \item if $(\sp f\setminus\{0\})\subset \Gamma$,
then $f(z)$ tends to a finite limit as $y\to\infty,\,y\in\Gamma'$,
uniformly in $x\in\R^p$ for all $\Gamma'=\overline{\Gamma}'
\subset\Int\widehat{\Gamma}\cup\{0\}$,
 \item if $(\sp f\setminus\{0\})\subset\Lambda+\Gamma$ with some
$\Lambda\in\sp f\cap(-\Gamma)\setminus\{0\}$, then  the function
$f(z)$ tends to $\infty$ as $y\to\infty,\,y\in\Gamma'$,  uniformly
in $x\in\R^p$ for all $\Gamma'=\overline{\Gamma}'
\subset\Int\widehat{\Gamma}\cup\{0\}$,
 \item if $(\sp f\setminus\{0\})\subset\Lambda+\Gamma$ with some
  $\Lambda\in(\overline{\sp f}\setminus{\sp f})\cap(-\Gamma)\setminus\{0\}$,
then  the function $f(z)$  takes every complex value on the set
$\Int T_{\Ga\cap\{|y|>q\}}$ for each $q<\infty$,
 \item if
$(\sp f\setminus\{0\})\subset\Lambda+\Gamma$  with some
 $\Lambda\in \overline{\sp f}\setminus\left((-\Gamma)\cup\Gamma\right)$,
  then the function $f(z)$ takes every complex value, except for at most one,
 on the set $\Int T_{\Ga\cap\{|y|>q\}}$ for each $q<\infty$,
 \item if $(\sp f\setminus\{0\})\not\subset\Lambda+\Gamma$
 for all $\Lambda\in \overline{\sp f}$ and
$\sp f\not\subset\Gamma$,  then the function $f(z)$  takes  every
complex value on the set $\Int T_{\Ga\cap\{|y|>q\}}$ for each
$q<\infty$.
\end{thm}

\medskip
{\bf Remark}. It is clear that we can replace $\sp
f\setminus\{0\}$ by $\sp f$ in cases 1 - 3. Therefore Theorem
\ref{thm2} gives,  in a sense,  a complete description of the
value distributions for our class of almost periodic functions.

\begin{proof}
 Case 1 was proved in \cite{FaUd}, case 2 was proved in \cite{gifa}.
Reduce case 3 to one-dimensional
one. Take $y^0\in\Int\widehat{\Gamma}$ such that $\langle
y^0,\,\lambda^{k}\rangle\neq\langle y^0,\,\lambda^{m}\rangle$ for
all $k\neq m$, and put $\varphi(w)=f(wy^0), w\in\CC$.

First, check that $\sp\varphi =\{ \langle y^0,\,\lambda \rangle :
\lambda\in\sp f\}$.  This is evident for finite exponential sums.
In the general case, take a sequence of Bochner-Feyer exponential
sums \footnote{For almost periodic functions on $\R$ see
\cite{lev}, Ch.VI, \S1, or \cite{C}, p.38-45; consideration in the
multidimensional case is similar.} $P_n(x)$, which approximates
$f(x)$ on $R^p$. Since $\sp P_n\subset\sp f$ and $P_n(uy^0)\to
\varphi(u)$ uniformly on $\R$, we see that
$\sp\varphi\subset\{\langle y^0,\, \lambda\rangle:\, \lambda\in\sp
f\}$. On the other hand, if $\lambda\in\sp f$, then
 $$
a_{\langle
y^0,\,\lambda\rangle}(0,\,P_n(y^0u))=a_{\lambda}(0,\,P_n)\to
a_\lambda(0,\,f)\neq 0\quad \hbox{as}\quad n\to\infty.
 $$
Therefore, $a_{\langle y^0,\,\lambda\rangle}(0,\,\varphi)\neq 0$
and ${\langle y^0,\,\lambda\rangle}\in\sp\varphi$.

Note that $\langle y^0,\,\lambda^{n}\rangle\to\langle
y^0,\,\Lambda\rangle$ as $\lambda^n\to\Lambda$, $\lambda^n\in\sp
f$. Also, since $y^0\in\Int\widehat{\Gamma}$ and
$\lambda-\Lambda\in\Gamma$ for all $\lambda\in\sp f$, we get
$\langle y^0,\,\lambda\rangle>\langle y^0,\,\Lambda\rangle$.
Therefore, $\inf\sp\varphi=\langle y^0,\,\Lambda\rangle$ and
$\langle y^0,\,\Lambda\rangle\not\in\sp\varphi$. From Theorem B,
i.3 it follows that $f(z)$ takes every complex value
 on the set $\{z=wy_0:\, \im w>q\}$ for each $q<\infty$.

Let us consider case 4. Let $b_0$ be the   coefficient of series
(\ref{Dir}) corresponding to the exponent $\lambda=0$. Then for
any $A\in\CC\setminus\{b_0\}$ each function $f(z)-A$ has the
spectrum $\sp f\cup\{0\}$.  Suppose that the support function
$H_{\sp f\cup\{0\}}(y)$ is linear on $(-\widehat{\Gamma})$. Then
it is not hard to prove (for example, see \cite{gifa}, Lemma 2)
that $\sp f\cup\{0\}\subset\Lambda'+\Gamma$ with some
$\Lambda'\in(-\Gamma)\cap(\overline{\sp f\cup\{0\}})$. But  this
is impossible in our case. Hence, the function $H_{\sp
f\cup\{0\}}(y)$ is nonlinear on $(-\widehat{\Gamma})$ Now
Corollary \ref{cor} yields that the function $f(z)-A$ has zeros on
$\Int T_{\Ga\cap\{|y|>q\}}$ for each $q<\infty$.

Let us consider case 5. Let $b_0$ be the same as in case 4. The
function $f(z)-b_0$ has the spectrum $\sp f\setminus\{0\}$. Note
that the support function $H_{\sp f\setminus\{0\}}(y)$ is
nonlinear on $(-\widehat{\Gamma})$. Hence Corollary \ref{cor}
implies that the function $f(z)-b_0$ has zeros on $\Int
T_{\Ga\cap\{|y|>q\}}$ for each $q<\infty$. Further, for any
$A\in\CC\setminus\{b_0\}$ the function $f(z)-A$ has the spectrum
$\sp f\cup\{0\}$. If the support function $H_{\sp f\cup\{0\}}(y)$
is linear on $(-\widehat{\Gamma})$, then $\sp
f\cup\{0\}\subset\Lambda'+\Gamma$ with some
$\Lambda'\in(-\Gamma)\cap(\overline{\sp f\cup\{0\}})$. The both
cases $\Lambda'=0$ and $\Lambda'\neq 0$ contradict to the
conditions of case 5.  Therefore the function $f(z)-A$ has zeros
on $\Int T_{\Ga\cap\{|y|>q\}}$ for each $q<\infty$.
\end{proof}

\end{document}